\documentclass[a4paper, 12pt]{article}

\usepackage[koi8-r,cp1251]{inputenc}
\usepackage{amsfonts,amssymb,amsmath, hyperref}
\usepackage[final]{epsfig}
\usepackage{graphicx}

\begin{document}

\begin{center}
\textbf{ HAUSDORFF OPERATORS  ON LEBESGUE SPACES WITH POSITIVE DEFINITE PERTURBATION MATRICES ARE NON-RIESZ}
\end{center}

\

\begin{center}
\textbf{A. R. Mirotin}
\end{center}

\begin{center}
amirotin@yandex.ru
\end{center}

\

\textsc{Abstract.} {\small We consider  generalized  Hausdorff operators with positive definite and permutable perturbation matrices on Lebesgue spaces and   prove that such  operators are not  Riesz operators provided they are non-zero.}

\


\

\section{Introduction and preliminaries}

The one-dimensional Hausdorff transformation
$$
(\mathcal{H}_1f)(x) =\int_\mathbb{R} f(xt)d\chi(t),\eqno(1)
$$
where $\chi$ is a measure on $\mathbb{R}$ with support $[0,1]$, was introduced by Hardy \cite[Section 11.18]{H} as a  continuous variable analog  of regular  Hausdorff transformations (or Hausdorff means) for series. Its modern $n$-dimensional generalization looks as follows:
$$
(\mathcal{H}f)(x) =\int_{\mathbb{R}^m} \Phi(u)f(A(u)x)du, \eqno(2)
$$
where $\Phi:\mathbb{R}^m\to \mathbb{C}$ is a locally integrable function,  $A(u)$ stands  for a family of non-singular $n\times n$-matrices, $x\in \mathbb{R}^n$, a column vector. See survey articles \cite{Ls}, \cite{CFW}  for historical remarks and the state of the art up to 2014.

To justify this definition the following approach may be suggested. Hardy \cite[Theorem 217]{H}
proved that (if  $\chi$ is a  probability measure) the transformation (1) gives rise to a regular generalized limit at infinity of the function $f$ in a sense that if $f$ is continuous  on $\mathbb{R},$ and $f(x) \to l$   then $\mathcal{H}_1f(x) \to l$ when $x\to \infty.$ Note that the map $x\mapsto xt$ ($t\ne 0$) is the general form of  automorphisms of the additive group $\mathbb{R}$.
This observation  leads to the  definition of a (generalized) Hausdorff  operator on a general group $G$ via the automorphisms of $G$ that was introduced and studied by the author in \cite{JMAA}, and \cite{AddJMAA}. For the additive group $\mathbb{R}^n$
this definition looks as follows.

\textbf{Definition 1.} Let  $(\Omega,\mu)$ be some $\sigma$-compact  topological space endowed with a positive regular  Borel measure $\mu,$ $\Phi$  a locally integrable function on $\Omega,$  and $(A(u))_{u\in \Omega}$  a $\mu$-measurable family of $n\times n$-matrices that are nonsingular
for $\mu$-almost every $u$ with $\Phi(u) \ne 0.$
  We define the \textit{Hausdorff  operator} with the kernel $\Phi$  by ($x\in\mathbb{R}^n$ is a column vector)
$$
(\mathcal{H}_{\Phi, A}f)(x) =\int_\Omega \Phi(u)f(A(u)x)d\mu(u).
$$

The general form of a Hausdorff  operator given by definition 1 (with an arbitrary measure space $(\Omega,\mu)$ instead of $\mathbb{R}^m$) gives us, for example, the opportunity to consider (in the case $\Omega=\mathbb{Z}^m$)
 discrete Hausdorff  operators  \cite{Forum}, \cite{faa}.

As was mentioned above Hardy proved that  the Hausdorff  operator (1) possesses some regularity property. For the operator given by the definition 1 the multidimensional version of his result is also true as the next proposition shows.

\textbf{Proposition 1.} \cite{faa} \textit{Let the conditions of definition} 1 \textit{are fulfilled.   In order that the transformation $\mathcal{H}_{\Phi, A}$  should be
regular, i.e. that $f$ is measurable and locally bounded on $\mathbb{R}^n,$ $f(x) \to l$ when $x\to \infty$   should imply $\mathcal{H}_{\Phi, A}f(x) \to l$, it is necessary and
sufficient that $\int_\Omega \Phi(u)d\mu(u)=1.$}

So, as for the classic transformation considered by Hardy
 the Hausdorff transformation  in the sense of the definition 1 gives rise to a new family (for various  $(\Omega,\mu)$, $\Phi$, and $A(u)$) of regular generalized limits at infinity for  functions of $n$ variables.

(For a different approach to justify the  definition (2) see \cite{LK}.)

The problem of compactness of  Hausdorff  operators was   posed by Liflyand \cite{L} (see also \cite{Ls}).
There is a conjecture that nontrivial Hausdorff operator in $L^p(\mathbb{R}^n)$ is non-compact.
 For  the case  $p=2$ and for commuting $A(u)$ this hypothesis was confirmed in \cite{Forum} (and for the diagonal $A(u)$ --- in \cite{JMAA}). Moreover, we  conjecture that every nontrivial Hausdorff operator in $L^p(\mathbb{R}^n)$ is non-Riesz.

Recall that a \textit{Riesz operator} $T$ is a bounded operator on some Banach space with spectral properties like those
of a compact operator; i.~e., $T$ is a non-invertible operator whose
nonzero spectrum  consists  of eigenvalues  of finite multiplicity with no limit points other then $0$. This is equivalent to the fact that $T-\lambda$ is Fredholm for every $\lambda \ne 0$ \cite{Ruston}. For  example, a sum of a quasinilpotent and compact operator is Riesz \cite[Theorem 3.29]{Dow}. Other interesting characterizations for Riesz operators one can also find in \cite{Dow}.

In this note we prove the  above mentioned conjecture for the case where the family $A(u)$ consists of  permutable and  positive (negative)  definite matrices.

\section{The main result}

We shall employ three lemmas to prove our main result.

\textbf{Lemma 1} \cite{JMAA} (cf. \cite[(11.18.4)]{H}, \cite{BM}). \textit{Let  $|\det A(u)|^{-1/p}\Phi(u)\in L^1(\Omega).$ Then the operator $\mathcal{H}_{\Phi, A}$ is bounded in   $L^p(\mathbb{R}^n)$ ($1\leq p\le \infty$) and}
$$
\|\mathcal{H}_{\Phi, A}\|\leq \int_\Omega |\Phi(u)||\det A(u)|^{-1/p}d\mu(u).
$$
This estimate is sharp  (see  theorem 1  in \cite{faa}).

\textbf{Lemma 2} \cite{faa} (cf. \cite{BM}). \textit{Under the conditions of Lemma} 1 \textit{the adjoint for the Hausdorff operator in $L^p(\mathbb{R}^n)$ has the form}
$$
(\mathcal{H}_{\Phi, A}^*f)(x) =\int_\Omega \overline{\Phi(v)}|\det A(v)|^{-1}f(A(v)^{-1}x)d\mu(v).
$$
\textit{Thus, the adjoint for a Hausdorff operator  is also  Hausdorff.}

\textbf{Lemma 3.} \textit{Let   $S$ be a boll in $\mathbb{R}^n$, $q\in [1,\infty)$, and $R_{q,S}$ denotes the restriction operator $L^q(\mathbb{R}^n)\to  L^q(S)$, $f\mapsto f|S$. If we as usual identify the dual of $L^q$ with $L^p$ ($1/p+1/q=1$), then the adjoint $R_{q,S}^*$ is the operator of natural embedding $L^p(S)\hookrightarrow L^p(\mathbb{R}^n)$}.

Proof. For $g\in L^p(S)$ let
$$
g^*(x)=\begin{cases}
g(x)\  \mathrm{for}\  x\in S,\\
0 \quad \quad \mathrm{for}\  x\in \mathbb{R}^n\setminus S.
\end{cases}
$$
Then the map $g\mapsto g^*$ is the natural embedding $L^p(S)\hookrightarrow L^p(\mathbb{R}^n)$.

By definition, the adjoint  $R_{q,S}^*: L^q(S)^*\to L^q(\mathbb{R}^n)^*$ acts according to  the rule
$$
(R_{q,S}^*\Lambda)(f)=\Lambda(R_{q,S}f)\quad \ (\Lambda\in L^q(S)^*, f\in L^q(\mathbb{R}^n)).
$$
If we (by the Riesz theorem) identify the dual of $L^q(S)$ with $L^p(S)$ via the formula $\Lambda\leftrightarrow g$, where
$$
\Lambda(h)=\int_S g(x)h(x)dx\quad \ (g\in L^p(S), h\in L^q(S)),
$$
and analogously for the the dual of $L^q(\mathbb{R}^n)$, then the definition of  $R_{q,S}^*$ takes the form
$$
\int_{\mathbb{R}^n}(R_{q,S}^*g)(x)f(x)dx=\int_S g(x)(f|S)(x)dx.
$$
But
$$
\int_S g(x)(f|S)(x)dx=\int_{\mathbb{R}^n}g^*(x)f(x)dx\quad (f\in L^q(\mathbb{R}^n)).
$$
 The right-hand side of the last formula is the linear functional from $L^q(\mathbb{R}^n)^*$.
 If we (again by the Riesz theorem) identify this functional with the function $g^*$, the result follows.$\Box$

\textbf{Theorem 1.}
\textit{Let $A(v)$ be a commuting family of real positive definite $n\times n$-matrices  ($v$ runs over the support of   $\Phi$), and $(\det A(v))^{-1/p}\Phi(v)\in L^1(\Omega).$ Then every nontrivial Hausdorff  operator $\mathcal{H}_{\Phi, A}$ in  $L^p(\mathbb{R}^n)$ ($1\leq p\le \infty$) is a non-Riesz  operator (and in particular it is non-compact).}

Proof. Assume the contrary. Since   $A(u)$ form a commuting family, there are an orthogonal $n\times n$-matrix $C$ and a family of
diagonal non-singular real matrices $A'(u)$ such that
$A'(u) = C^{-1}A(u)C$ for $u\in \Omega.$ Consider the bounded and invertible operator
 $\widehat{C}f(x):=f(Cx)$ in $L^p(\mathbb{R}^n).$ Because of the equality $\widehat{C}\mathcal{H}_{\Phi, A}\widehat{C}^{-1}=\mathcal{H}_{\Phi, A'},$ operator $\mathcal{H}:=\mathcal{H}_{\Phi, A'}$ is Riesz and nontrivial, too.

Note that each open hyperoctant in   $\mathbb{R}^n$  is $A(u)$-invariant. Chose such an open
 $n$-hyperoctant  $U$ that   $\mathcal{K}:=\mathcal{H}|L^p(U)\ne 0.$ Then $L^p(U)$  is a closed $\mathcal{K}$-invariant subspace of $L^p(\mathbb{R}^n)$ and $\mathcal{K}$ is a nontrivial Riesz operator in  $L^p(U)$ by \cite[p. 80, Theorem 3.21]{Dow}.

Let $1\leq p<\infty.$ To get a contradiction, we shall use the modified $n$-dimensional  Mellin transform for the $n$-hyperoctant $U$ in the form

$$
(\mathcal{M}f)(s):=\frac{1}{(2\pi)^{n/2}}\int_{U}|x|^{-\frac{1}{q}+is}f(x)dx,\quad s\in \mathbb{R}^n
$$
  Here and below we assume that $|x|^{-\frac{1}{q}+ is}$ $:= \prod_{j=1}^n |x_j|^{-\frac{1}{q}+ is_j}$ where $|x_j|^{-\frac{1}{q}+ is_j}:=$ $\exp((-\frac{1}{q}+ is_j)\log |x_j|)$.
 The map $\mathcal{M}$ is a bounded operator between  $L^p(U)$ and $L^q(\mathbb{R}^n)$ for $1\leq p\leq 2$  ($1/p+1/q=1).$ It can be easily obtained from the Hausdorff--Young inequality for the $n$-dimensional  Fourier transform by using the exponential change of variables
 (see \cite{BPT}). Let $f\in L^p(U).$ First assume that   $|y|^{-1/q}f(y)\in L^1(U).$ Then
  as in the proof of  theorem 1 from \cite{Forum}, using the Fubini--Tonelli's theorem and integrating by substitution     $x=A(u)'^{-1}y,$ yield the following
  $$
  (\mathcal{MK}f)(s)=
\varphi(s)(\mathcal{M}f)(s)\  (s\in \mathbb{R}^n),
$$
where the function $\varphi$ ("the symbol of the the Hausdorff operator" \cite{Forum}) is bounded and continuous on    $\mathbb{R}^n.$

Thus,
$$
\mathcal{MK}f=\varphi \mathcal{M}f. \eqno(3)
$$
By continuity the last equality is valid for all  $f\in L^p(U).$

  Let $1\leq p\leq 2.$ There exists a constant $c> 0,$ such that the set    $\{s\in \mathbb{R}^n: |\varphi(s)|>c\}$ contains an open ball  $S.$ Formula  (3) implies that
$$
M_{\psi}R_{q,S}\mathcal{M}\mathcal{K}=R_{q,S}\mathcal{M},
$$
  where $\psi=(1/\varphi)|S,$  $M_{\psi}$ denotes the operator of multiplication by  $\psi$, and $R_{q,S}: L^q(\mathbb{R}^n)\to L^q(S),$ $f\mapsto f|S$ --- the restriction  operator. Let $T=R_{q,S}\mathcal{M}.$ Passing to the conjugates  gives
 $$
 \mathcal{K}^*T^*M_{\psi}^*=T^*.
 $$
    By \cite[Theorem 1]{FNR3} this implies that the operator $T^*=\mathcal{M}^*R_{q,S}^*$  has finite rank. By Lemma 3 $R_{q,S}^*$ is the operator of natural embedding $L^p(S)\hookrightarrow L^p(\mathbb{R}^n)$.

    For  $g\in L^p(\mathbb{R}^n)$ consider the operator
$$
(\mathcal{M}'g)(x):=\frac{1}{(2\pi)^{n/2}}\int_{\mathbb{R}^n}|x|^{-\frac{1}{q}+is}g(s)ds,\  x\in U.
$$
This  is a bounded operator taking  $L^p(\mathbb{R}^n)$ into  $L^q(U).$ Indeed,
since
$$
|x|^{-\frac{1}{q}+ is}=\prod_{j=1}^n |x_j|^{-\frac{1}{q}}\exp(is_j\log |x_j|),
$$
we have
$$
(\mathcal{M}'g)(x)=|x|^{-\frac{1}{q}}\frac{1}{(2\pi)^{n/2}}\int_{\mathbb{R}^n}\exp(is\cdot \log|x|)g(s)ds,\  x\in U,
$$
where $|x|:=|x_1|\dots |x_n|$, $\log |x|:=(\log |x_1|,\dots,\log |x_n|)$, and the dot denotes the inner product in $\mathbb{R}^n$. Thus, we can express the function  $\mathcal{M}'g$ via the Fourier transform $\widehat{g}$ of $g$ as follows: $(\mathcal{M}'g)(x)=|x|^{-1/q}\widehat{g}(-\log |x|)$,  ($x\in U$) and therefore
$$
\|\mathcal{M}'g\|_{L^q(U)}=\left(\int_U |x|^{-1}|\widehat{g}(-\log |x|)|^q dx\right)^{1/q}.
$$
Putting here $y_j:=-\log|x_j|$ ($j=1,\dots,n$) and taking into account that the Jacobian of this transformation is
$$
\frac{\partial(x_1,\dots,x_n)}{\partial(y_1,\dots,y_n)}=\det \mathrm{diag}(e^{-y_1},\dots,e^{-y_n})=\exp\left(-\sum_{j=1}^n y_j \right),
$$
we get by  the Hausdorff--Young inequality that
$$
\|\mathcal{M}'g\|_{L^q(U)}=\|\widehat{g}\|_{L^q(\mathbb{R}^n)}\le \|g\|_{L^p(\mathbb{R}^n)}.
$$

If  $f\in L^p(U)$, and  $f(x)|x|^{-1/q}\in L^1(U)$, $g\in  L^p(\mathbb{R}^n)\cap L^1(\mathbb{R}^n)$ the Fubini--Tonelli’s theorem implies
$$
\int_{\mathbb{R}^n}(\mathcal{M}f)(s)g(s)ds=\int_{U}f(x)(\mathcal{M}'g)(x)dx.
$$
Since the bilinear form $(\varphi,\psi)\mapsto\int \varphi \psi d\mu$ is continuous on $L^p(\mu)\times L^q(\mu)$, the last equality is valid for all   $f\in L^p(U)$, $g\in  L^p(\mathbb{R}^n)$ by continuity. So, $\mathcal{M}'=\mathcal{M}^*$.

  It was shown above that the restriction of the operator  $\mathcal{M}^*$  to  $L^p(S)$ has finite rank. Since $\mathcal{M}^*$ can be easily reduced to the Fourier transform, this is contrary to the Paley--Wiener theorem on the Fourier image of the space  $L^2(S)$ ($L^2(S)\subseteq L^p(S)$) (see, i.~g., \cite[Theorem  III.4.9]{SW}).

  Finally, if $2<p\le \infty$ one can use duality arguments. Indeed, by lemma 2 the adjoint operator $\mathcal{H}_{\Phi, A'}^*$ (as an operator in $L^q(\mathbb{R}^n)$) is also of Hausdorff type. More precisely, it equals to  $\mathcal{H}_{\Psi, B},$ where $B(u)= A(u)'^{-1}= \mathrm{diag}(1/a_1(u),\dots,1/a_n(u))$,  $\Psi(u)=\Phi(u)|\det A(u)'^{-1}|=$ $\Phi(u)/a(u)$.  It is easy to verify that  $\mathcal{H}_{\Psi, B}$ satisfies all the  conditions of theorem 1  (with  $q,$ $\Psi$ and $B$ in place of $p,$ $\Phi$ and $A$ respectively).
   Since $1\le q<2$, the operator  $\mathcal{H}_{\Psi, B}$ is not a Riesz operator in $L^q(\mathbb{R}^n),$ and so is $\mathcal{H}_{\Phi, A}$,  because   $T$ is a Riesz operator if only if its conjugate $T^*$ is a
Riesz operator  \cite[p. 81, Theorem 3.22]{Dow}.$\Box$

\section{Corollaries and examples}

For the next corollary we need the following

\textbf{Lemma 4}. \textit{Let $J:X\to X$ be a linear isometry of a Banach space $X$. A bounded  operator $T$ on $X$ which commutes with $J$ is a Riesz operator if and only if such is $JT$.}

Proof. We use the fact that an operator $T$ is a Riesz operator if and only if it is
asymptotically quasi-compact \cite {Ruston} (see also \cite[Theorem 3.12]{Dow}). This means that
$$
\lim\limits_{n\to\infty}\left(\inf\limits_{C\in \mathcal{K}(X)}\|T^n-C\|^{1/n}\right)=0,
$$
where $\mathcal{K}(X)$ denotes the ideal of compact operators in $X$  (Ruston condition). Since $(UT)^n=U^nT^n$ and
$$
\inf\limits_{C\in \mathcal{K}(X)}\|(UT)^n-C\|^{1/n}=\inf\limits_{C\in \mathcal{K}(X)}\|T^n-U^{-n}C\|^{1/n}=\inf\limits_{C'\in\mathcal{ K}(X)}\|T^n-C'\|^{1/n},
$$
the result follows.$\Box$

\textbf{Corollary 1}. \textit{Let $A(v)$ be a commuting family of real negative definite $n\times n$-matrices  ($v$ runs over the support of   $\Phi$), and $(\det A(v))^{-1/p}\Phi(v)\in L^1(\Omega).$ Then every nontrivial Hausdorff  operator $\mathcal{H}_{\Phi, A}$ in  $L^p(\mathbb{R}^n)$ ($1\leq p\le\infty$) is non-Riesz  (and in particular it is non-compact).}

Proof. Let $Jf(x):=f(-x)$.  Since $-A(v)$ form a commuting family of real positive definite $n\times n$-matrices, and $\mathcal{H}_{\Phi, A}=J\mathcal{H}_{\Phi, -A}$, this corollary follows from lemma 4 and theorem 1. $\Box$

\textbf{Corollary 2.} \textit{Under the conditions of theorem 1 or corollary 1 Hausdorff  operator $H_{\Phi, A}$ is not the sum of the quasinilpotent and compact operators.}

Indeed, as was mentioned in the introduction, the sum of the quasinilpotent and compact operators is a Riesz operator.

\textbf{Corollary 3.} \textit{Let $n=1$, $\phi:\Omega\to \mathbb{C}$ and let $a(v)$ be a  real and positive (negative) function on $\Omega$  ($v$ runs over the support of   $\phi$), and $|a(v)|^{-1/p}\phi(v)\in L^1(\Omega).$ Then every nontrivial Hausdorff  operator
$$
\mathcal{H}_{\phi, a}f(x)=\int_{\Omega} \phi(u)f(a(u)x)d\mu(u)\  (x\in \mathbb{R})
$$
in  $L^p(\mathbb{R})$ ($1\leq p\le\infty$) is a non-Riesz  operator (and in particular it is non-compact).}

\textbf{Example 1}. Let $t^{-1/q}\psi(t)\in L^1(0,\infty).$
Then by corollary 3 the operator
$$
\mathcal{H}_\psi f(x)=\int_0^\infty\frac{\psi(t)}{t}f\left(\frac{x}{t}\right)dt
$$
is a non-Riesz  operator in  $L^p(\mathbb{R})$ ($1\leq p\le\infty$) provided it is non-zero.

\textbf{Example 2}. Let $(t_1t_2)^{-1/p}\psi_2(t_1,t_2)\in L^1(\mathbb{R}_+^2).$
Then by theorem 1 the operator
$$
\mathcal{H}_{\psi_2} f(x_1,x_2)=\frac{1}{x_1x_2}\int_0^\infty\!\int_0^\infty\psi_2\left(\frac{t_1}{x_1}, \frac{t_2}{x_2}\right) f(t_1,t_2)dt_1dt_2
$$
is a non-Riesz  operator in  $L^p(\mathbb{R}_+^2)$ ($1\leq p\le\infty$) provided it is non-zero.


\begin{thebibliography}{99}



\bibitem{H}
G. H. Hardy, Divergent series, Clarendon Press, Oxford, 1949.


\bibitem{Ls}
E. Liflyand, Hausdorff operators on Hardy spaces, Eurasian Math.  J. ,
no. 4,  101 -- 141  (2013)


\bibitem{CFW}
J. Chen, D. Fan, S. Wang, Hausdorff operators on Euclidean space (a survey article), Appl. Math. J. Chinese Univ. Ser. B (4), 28,  548--564 (2014)


\bibitem{JMAA}
A. R. Mirotin,Boundedness of Hausdorff operators on
Hardy spaces $H^1$ over locally compact groups,  J. Math. Anal. Appl., \textbf{473} (2019), 519 -- 533. DOI  10.1016/j.jmaa.2018.12.065. Preprint arXiv:1808.08257v2 [math.FA] 1 Sep 2018.

\bibitem{AddJMAA}
A. R. Mirotin, Addendum to "Boundedness of Hausdorff operators on
Hardy spaces $H^1$ over locally compact groups", J. Math. Anal. Appl., vol. 479, No. 1,  872 -- 874  (2019).



\bibitem{LK}
E. Liflyand, A. Karapetyants,
Defining Hausdorff operators on Euclidean spaces,
Mathematical Methods in the Applied Sciences, 2020, DOI: 10.1002/mma.6448.


\bibitem{Forum}
A. R. Mirotin,  The structure of normal  Hausdorff operators  on Lebesgue spaces, Forum Math.,  2020 - V. 32, No 1 - P. 111-119. https://doi.org/10.1515/forum-2019-0097.



\bibitem{faa}
A. R. Mirotin, On the Structure of Normal Hausdorff Operators
on Lebesgue Spaces, Functional Analysis and Its Applications, 2019, Vol. 53, No. 4, pp. 261--269.




\bibitem{L}
E. Liflyand,  Open problems on Hausdorff operators, In: Complex Analysis and Potential Theory, Proc. Conf. Satellite to ICM 2006, Gebze, Turkey, 8-14 Sept. 2006; Eds. T. Aliyev Azeroglu and P.M. Tamrazov; World Sci., 280--285 (2007)

\bibitem{Ruston}
A. F. Ruston, Operators with Fredholm theory, J. London Math. Soc. 29 (1954),
318--326.

\bibitem{Dow}
H. R. Dowson, Spectral Theory of Linear Operators, Academic Press Inc., London, 1978.



\bibitem{BM}
G. Brown and F. M\'{o}ricz, Multivariate Hausdorff operators on the spaces $L^p(\mathbb{R}^n),$
J. Math. Anal. Appl., 271, 443--454  (2002)




\bibitem{BPT}
Yu.A. Brychkov, H.-J. Glaeske, A.P. Prudnikov, and Vu Kim Tuan, Multidimentional
Integral Transformations. Gordon and Breach, New York - Philadelphia - London - Paris -
Montreux - Tokyo - Melbourne - Singapore, 1992.



\bibitem{FNR3}
C. K. Fong,
E. A. Nordgren,
M. Radjabalipour, H. Radjavi, P. Rosenthal,
Extensions of Lomonosov's invariant subspace theorem. Acta Sci. Math. (Szeged), 41,  55 -- 62 (1979)



\bibitem{Hilb}
A. R. Mirotin, Hilbert transform in context
of locally compact abelian groups, Int. J. Pure and Appl.
Math., 51,  463 -- 474 ( 2009)



\bibitem{SW}
E. M. Stein, G. Weiss, Introduction to Fourier Analysis on Eucledean Spaces, Prinston Univercity Press, Prinston, New Jersey, 1971.








\end{thebibliography}
\end{document}